\documentclass{amsart} 
\usepackage{fullpage}

\usepackage[alphabetic]{amsrefs}

\newcommand{\spacesymbol}[1]{\mathbb{#1}}
\newcommand{\R}{\spacesymbol{R}}

\newcommand{\abs}[1]{\lvert #1 \rvert}
\newcommand{\norm}[1]{\lVert #1 \rVert}

\newcommand{\m}[2]{\left\langle #1, #2 \right\rangle}

\DeclareMathOperator\tr{tr}

\DeclareMathOperator\ad{ad}

\DeclareMathOperator\ric{Rc}

\DeclareMathOperator\GL{GL}

\DeclareMathOperator\Der{Der}

\newtheorem{definition}{Definition}
\newtheorem{prop}{Proposition}
\newtheorem{theorem}{Theorem}
\newtheorem{corollary}{Corollary}
\newtheorem{lemma}{Lemma}

\theoremstyle{remark}
\newtheorem*{remark}{Remark}
\newtheorem*{acknowledgements}{Acknowledgements}

\begin{document} 

\title{Bach Flow of Simply Connected Nilmanifolds}
\author{Adam Thompson}
\begin{abstract}
    The Bach flow is a fourth order geometric flow defined on four manifolds.  For a compact manifold, it is a conformally modified gradient flow for the \(L^2\)-norm of the Weyl curvature.  In this paper we study the Bach flow on four-dimensional simply connected nilmanifolds whose Lie algebra is indecomposable.  We show that the Bach flow beginning at an arbitrary left invariant metric exists for all positive times and after rescaling converges in the pointed Cheeger-Gromov sense to an expanding Bach soliton which is non-gradient.  Combining our results with previous results of Helliwell gives a complete description of the Bach flow on simply connected nilmanifolds.
\end{abstract}

\maketitle
\section{Introduction}

The success of the Ricci flow and other second order geometric flows has led to a desire to better understand higher order geometric flows. In particular, the study of fourth order flows is a natural progression from second order flows.  Fourth order flows can be more difficult to since we no longer have access to maximum principles. An important class of fourth order flows are the gradient flows of quadratic curvature functionals.

In \cite{BahuaudHelliwellSTEFSHOGF} Bahuad and Helliwell introduced a family of higher order flows called the \emph{ambient obstruction flow} which involves the ambient obstruction tensor of Fefferman and Graham introduced in~\cite{FeffermanGrahamCI}.  In dimension 4, the ambient obstruction tensor coincides with the \emph{Bach tensor} and the corresponding flow is the \emph{Bach flow}: \begin{equation}\label{Eq:BachFLow}
\frac{\partial}{\partial t}g = \operatorname{B}(g)+\frac{\Delta s}{12}g, \qquad g(0)=g_0.
\end{equation} 
Here \(\operatorname{B}(g)\) is the Bach tensor defined in \eqref{Eq:BachTensor} below, \(s=s(g)\) is the scalar curvature and \(\Delta\) is the Laplace-Beltrami operator.  This flow corresponds to a conformally modified gradient flow of the \(L^2\)-norm of the Weyl curvature, and so is fourth order (see the introduction in \cite{StreetsTLTBOFOCF}). While this might not be the first choice of a Bach flow due to the \(\frac{1}{12}(\Delta s)g\) term, this is needed to counteract the invariance of the Bach tensor under the conformal group \cite[Section 1.2]{LopezAOF}.  In any case, if a solution \(g(t)\) to \eqref{Eq:BachFLow} has constant scalar curvature at each point in time, as is the case for homogeneous solutions, then \(\Delta s\equiv 0\). Short time existence and uniqueness for the flow on compact manifolds was proved by Bahaud and Helliwell in \cite{BahuaudHelliwellSTEFSHOGF,BahuaudHelliwellUFSHOGF}.

The Bach flow was studied by Helliwell in \cite{HelliwellBFOHP} on homogeneous manifolds of the form \(M^4=S^1\times N^3\) where \(N^3\) is a closed, homogeneous 3-manifold and \(M\) is equipped with the product metric.  By lifting to the universal cover, this reduces to the study of the Bach flow on 9 simply connected homogeneous spaces of the form \(\widetilde{M}=\R\times \widetilde{N}\).  In \cite{GriffinGAOSOHM}, Griffin uses the computations of the Bach tensor from \cite{HelliwellBFOHP} and a theorem of Petersen and Wylie \cite[Theorem 3.6]{PetersenWylieROHGSMARE} to study gradient Bach solitons on the product manifolds from Helliwell's paper.

The Bach tensor of a left-invariant metric is determined by its value at the identity, and therefore can be computed explicitly using metric coefficients and structure constants of the of the Lie algebra. In practice, however, this is difficult due to the fourth order nature of the tensor.  For Riemannian products this formula greatly simplifies \cite[Section 2.1]{HelliwellBFOHP}.  In order to make the problem more tractable in the non-product case we have chosen to study the flow on Lie group whose Lie algebra has a large number of structure constants which vanish.

In this article we study the Bach flow the simply connected nilpotent Lie group, \(N^4\), whose Lie algebra is the unique four-dimensional indecomposable nilpotent Lie algebra.  By indecomposable, we mean that the Lie algebra cannot be written as a direct product of two lower dimensional Lie algebras.  Our main result is the following theorem.

\begin{theorem}\label{Thm:MainThm}

Let \(N^4\) be the simply connected nilpotent Lie group with Lie algebra \(\operatorname{Lie}(N)=\mathfrak{n}_4\), the unique four dimensional, indecomposoable, nilpotent Lie algebra and let \(g\) be a left invariant metric on \(N\).
Then, there exists a unique solution \((g(t))_{t\in (-\varepsilon,T)}\) to the Bach flow, defined on a maximal interval of existence \((-\varepsilon,T)\), which is \(N^4\)-invariant.  This satisfies
\begin{enumerate}
    \item \(T=\infty\) and the metrics \(g(t)\) converges in the pointed Cheeger-Gromov sense to the Euclidean metric in \(\R^4\) as \(t\to \infty\).
    \item If \(s_t=s(g(t))\) denotes the scalar curvature of \(g(t)\) then the rescaled metrics \[\tilde{g}(t)=|s_t|g(t)\] converge in the pointed Cheeger-Gromov sense to an expanding Bach soliton \(g_\infty\) on \(N^4\) as \(t\to\infty\). 
\end{enumerate}

\end{theorem}

A \emph{Bach soliton} is a solution to \eqref{Eq:BachFLow} which is self similar.  That is, the solution \(g(t)\) to \eqref{Eq:BachFLow} beginning at a soliton \(g\) can be written \(g(t)=c(t)\varphi(t)^*g\) for some positive function \(c\) and time dependent diffeomorphisms \(\varphi(t)\).  We say a soliton is expanding, steady or shrinking if the scaling constant is increasing, constant or decreasing respectively.

The Bach soliton \(g_\infty\) is a non-trivial algebraic Bach soliton (see \cite{LauretRSHN,LauretGFATSOHS} for literature on algebraic Ricci solitons).  That is, if \(B:\mathfrak{n}_4\to \mathfrak{n}_4\) is the endomorphism defined by \(\operatorname{B}(g_\infty)_e = g(B\cdot,\cdot)\) then, \[B = \lambda I +D,\qquad \lambda\in \R,\] where \(D\in \Der(\mathfrak{n}_4)\) is a non-zero derivation. To prove Theorem~\ref{Thm:MainThm} we use the bracket flow technique introduced by Lauret to study the Ricci flow in homogeneous manifolds \cite{LauretRFOHM} (Lauret extends the bracket flow technique to more general geometric structures in \cite{LauretGFATSOHS}).  The bracket Bach flow is an ODE in the variety of Nilpotent Lie brackets \(\mathcal{N}^4\subset \Lambda^2(\R^4)^*\otimes \R^4\) which is equivalent to the Bach flow.  Algebraic solitons correspond to solutions of the bracket flow which evolve only by scaling. These are easier to detect than general solitons since we have removed the action of the diffeomorphism group.  Our analysis follows Lauret's analysis of the Ricci flow on simply-connected nilmanifolds in \cite{LauretTRFFSCN}.  

We also remark that the soliton \(g_\infty\) is not of gradient type due to a Theorem of Petersen and Wylie \cite[Theorem 3.6]{PetersenWylieROHGSMARE}.  To our knowledge, this is the first example of a non-gradient Bach soliton.

The paper is structured as follows.  In \S \ref{Sec:Background} we introduce the Bach tensor and its basic properties.  In \S\ref{4DNilpotentLieAlgebras} review some group actions on the variety of nilpotent Lie algebras and recall the classification of four-dimensional Lie algebras. In \S\ref{Sec:LieGrp} we review the structure of four-dimensional simply connected nilmanifolds.  In particular, we show that up to isometry, left-invariant metrics on \(N^4\) depend on 3 real parameters.  By suitably gauging the Bach flow (Proposition~\ref{GuagedBracketFlow}), we obtain an ODE in three variables, \(a,b,c\), which is equivalent to the bracket Bach flow the simply connected Lie group \(N^4\).  By studying the long time behaviour of the quantities \(a,b,c\) we are able to understand the long time behaviour of the Bach flow.
Self similar solutions of the Bach flow are studied in \S\ref{Sec:Solitons}.  In terms of the variables \(a,b,c\) the soliton condition becomes a set of algebraic equations.  Finally, we considered the normalised flow in \S\ref{Sec:NormalisedFLows}.

\section{Background on The Bach Flow}\label{Sec:Background}

Let \(M^4\) be a four dimensional Riemannian manifold.  The \emph{Bach tensor} is the tensor given in local coordinates by \begin{equation}\label{Eq:BachTensor}
\operatorname{B}(g)_{ij} = \nabla^k\nabla^l W_{kijl}+\frac{1}{2}R^{kl}W_{ikjl}.
\end{equation}
Here, \(W\) is the Weyl tensor, \(R^{kl}=g^{ka}g^{lb}R_{ab}=g^{ka}g^{lb}(\ric)_{ab}\) are the components of the Ricci tensor with both indices raised.

 The Bach tensor is conformally invariant (\(\operatorname{B}(\rho g) = \rho^{-1}\operatorname{B}(g)\) for any positive function \(\rho\)) and arises in the study of manifolds which are conformally Einstein \cite[Section 6]{KuhnelRademacherCTOPRM}.  In particular, Theorem 6.6 in \cite{KuhnelRademacherCTOPRM} says \(\operatorname{B}(g)\equiv 0\) is a necessary condition for \((M^4,g)\) to be conformal to an Einstein manifold.  The converse holds if \(g\) is conformal to a metric with harmonic Weyl tensor \cite[Corollary 6.8]{KuhnelRademacherCTOPRM}.  The Bach tensor is trace and divergence free \cite{GriffinGAOSOHM}.

On a compact manifold the Bach tensor is \(-1/4\) times the gradient of the Riemannian functional \begin{equation*}
g \longmapsto \int_M |W|_g^2\, dv(g)
\end{equation*}
\cite[4.76]{BesseEM} (note that the definition given here is only agrees with the definition in \cite{BesseEM} up to a factor of \(-4\)).

The Bach flow has been studied on homogeneous \(4\)-manifolds whose universal cover has a product structure (that is, \((\widetilde{M}^4,g)=(N_1\times N_2, g_1\oplus g_2)\) where \((N_i,g_i)\) are simply connected homogeneous manifolds of lower dimension which admit compact quotients) by Helliwell in \cite{HelliwellBFOHP}. When \(\dim N_2=3\) we must have that \(N_2\) is a unimodular Lie group and that \(N_1=\R\).  In this case, Helliwell uses Milnor frames to diagonalise the Ricci and Bach tensor and then studies the resulting ODE's.  In particular, this includes the study of the Bach flow on \((\R^4,\bar{g})\) and \((\R \times H^3,\bar{g}\oplus g)\) where \(\bar{g}\) denotes the Euclidean metric on \(\R^n\) and \(g\) is a left invariant metric on the three-dimensional Heisenberg group \(H^3\).  For \(\dim N_1=\dim N_2=2\) the only possibilities are \(M(c_1)^2\times M(c_2)\) where \(M(c)^2\) denotes the simply connected space form with constant curvature \(c\in \R\).  In this case the flow remains a product of two space forms and the flow is equivalent to a system of ODE's for the curvature of the factors.

\section{Four-Dimensional Nilpotent Lie algebras}\label{4DNilpotentLieAlgebras}

Any bilinear skew-symmetric map \(\mu:\R^4\times \R^4\to \R^4\) which satisfies the \emph{Jacobi identity} defines a Lie algebra structure structure on \(\R^4\).  
There is a natural `change of basis' action of \(\GL_4(\R)\) on \(\Lambda^2(\R^4)^*\otimes \R^4\), the vector space of bilinear skew-symmetric maps, given by \begin{equation}\label{Eq:GLRep}
h\cdot \mu = h\mu(h^{-1}\cdot,h^{-1}\cdot), \qquad \forall h\in \GL_4(\R). 
 \end{equation}
Two Lie brackets \(\mu_1\) and \(\mu_2\) define isomorphic Lie algebras if and only if there is a \(h\in \GL_4(\R)\) such that \(h\cdot \mu_1=\mu_2\).

Moreover, any inner product \(\m{\cdot}{\cdot}\) on \(\R^4\) induces an inner product on \(\Lambda^2(\R^4)^*\otimes \R^4\) \cite{LauretTRFFSCN}: \begin{equation}\label{Eq:InnerProduct}
    \m{\mu}{\lambda} = \sum_{i,j=1}^4 \m{\mu(e_i,e_j)}{\lambda(e_i,e_j)}, \, \forall\mu,\lambda\in\Lambda^2(\R^n)^*\otimes \R^n
\end{equation} where \(\{e_i\}_{i=1}^4\) is any orthonormal basis.

If the map \((\ad_\mu x)\), defined by \((\ad_\mu x)y:=\mu(x,y)\) for all \(y\in \R^4\), is a nilpotent endomorphism for all \(x\in \R^4\), we say that \(\mu\) is \emph{nilpotent}.
A bracket \(\mu\) which satisfies the Jacobi identity and is nilpotent defines a nilpotent Lie algebra \cite[Theorem 3.5.4]{VaradarajanLGLAATR}.  In this way, the set \begin{equation}\label{EqVarietyOfNilLieAlgs}
\mathcal{N}_4 = \{\mu\in \Lambda^2(\R^4)^*\otimes \R^4: \mu\text{ is nilpotent and satisfies the Jacobi Identity.}\}
\end{equation} parametrises structure constants of \(4\)-dimensional, nilpotent Lie algebras. 

For any element \(\mu\in \mathcal{N}_4\), \((\R^4,\mu)\) is isomorphic to one of the three following Lie algebras \cite[Section 5.5]{deGraafCOSLA}: \begin{enumerate}
\item \(\R^4\): with the abelien Lie bracket \(\mu(x,y)=0\) for all \(x,y \in \R^4\).
\item \(\R\oplus \mathfrak{h}_3\): This is the product of the 3 dimensional Heisenberg algebra with \(\R\).  The non-trivial bracket relations are \(\mu(e_1,e_2)=-\mu(e_2,e_1)=e_3\).
\item \(\mathfrak{n}_4\): This is not a product Lie algebra.  The non-trivial bracket relations are \[\mu(e_1,e_2)=-\mu(e_2,e_1)=e_3,\quad \mu(e_1,e_3)=-\mu(e_3,e_1)=e_4.\]
\end{enumerate}
    Note that 1. and 2. are decomposable Lie algebras.

\section{Four-Dimensional Nilpotent Lie Groups}\label{Sec:LieGrp}

In this paper, we are concerned with simply connected Lie groups, \(N^4\), such that \(\operatorname{Lie}(N)\simeq \mathfrak{n}_4\).  Any such simply connected Lie group can be identified with \(\R^4\) under the operation \[x\cdot_\mu y := x+y+\frac{1}{2}\mu (x,y) +\frac{1}{12}\mu(\mu(x,y),y)-\frac{1}{12}\mu(\mu(x,y),x) \qquad \forall x,y\in \R^n.\] See \cite[Section 3]{VaradarajanLGLAATR}.

Let \(\m{\cdot}{\cdot}\) be the standard inner product on \(\R^n\).  Then any other inner product on \(\R^n\) can be written as \[(\cdot,\cdot)= \m{h\cdot}{h\cdot} \] for some \(h\in \GL_n(\R)\).  If \((\cdot,\cdot)\) is an inner product on \(\R^n\), we denote by \(g_{\mu,(\cdot,\cdot)}\) the left-invariant metric on \((\R^n,\cdot_\mu)\) which agrees with \((\cdot,\cdot)\) on \(T_0\R^n\) \cite{LauretTRFFSCN}.  When \((\cdot,\cdot)\) is the standard inner product on \(\R^n\), we will write \(g_\mu \) instead of \(g_{\mu,(\cdot,\cdot)}\). We then have the following proposition.

\begin{prop}[\cite{LauretTRFFSCN}, Theorem 4.1]\label{Metrics}
For \(\mu,\lambda\in \mathcal{N}_n\) and an inner product \((\cdot,\cdot)=\m{h\cdot}{h\cdot}\) on \(\R^n\), the metrics \(g_\mu\) and \(g_{\lambda,(\cdot,\cdot)}\) are isometric if and only if \(\lambda = h\cdot \mu\).  In particular, \(g_\mu\) and \(g_\lambda\) are isometric if and only if \(\lambda\in O(n)\cdot \mu\).
\end{prop}

Note that by the previous proposition, up to isometry we may assume that a left invariant metric agrees with the standard inner product on \(T_0\R^n\).

 Suppose now that \(\mu\in \mathcal{N}_4\) is a Lie bracket such that \((\R^4,\mu)\simeq \mathfrak{n}_4\) and \(\m{\cdot}{\cdot}\) is the standard inner product on \(\R^4\).  Proposition~\ref{Prop:Isom} below gives a convenient basis of \(\R^4\) for us to work in.

 \begin{prop}[\cite{VanThuongMO4DULG},Theorem 4.2]\label{Prop:Isom}
Let \((\R^4, g_{\tilde{\mu},(\cdot,\cdot)})\) be a four dimensional, simply connected Nilpotent Lie group equipped with a left-invariant metric \(g_{\tilde{\mu},(\cdot,\cdot)}\) such that \((\R^4,\tilde{\mu})\simeq \mathfrak{n}_4\).  Then \((\R^4, g_{\tilde{\mu},(\cdot,\cdot)})\) is isometric to \((\R^4,g_\mu)\) where the bracket \(\mu\) has the structure constants \begin{equation}\label{StructureConstants}
 \mu (e_1,e_2)=ae_3+be_4, \qquad \mu(e_1,e_3)=ce_4, \qquad a,b,c \in \R, a,c>0,
\end{equation} and \(\{e_i\}\) is the standard basis of \(\R^4\).
\end{prop}
 If \(\mu\in \mathcal{N}_4\) has the structure coefficients \eqref{StructureConstants} with respect to the standard basis of \(\R^4\) then we will write \(\mu=\mu_{a,b,c}\). 

Let us define
\begin{equation}
\mathcal{O} = \{\mu\in \mathcal{N}_4: \mu=\mu_{a,b,c}\text{ with respect to the standard basis of }\R^4\}.
\end{equation}
 It follows from the above discussion that if \((\R^n,g_\mu)\) is a simply connected Nilpotent Lie group with left-invariant metric \(g_\mu\), then up to isometry we may assume that \(\mu\in \mathcal{O}\).  Clearly we can identify \(\mathcal{O}\) with the set \[\{(a,b,c)\in \R^3: 0< a,c\}.\]

 We will use this in \S\ref{BachFlowSection} in order to reduce our study of the bracket flow, to be introduced later on.  Of course, there is no guarantee that a solution to the bracket flow \eqref{BracketFLow} will remain in \(\mathcal{O}\). However, we will show in \S\ref{BachFlowSection} that we can gauge our flow so that the solution does remain within \(\mathcal{O}\). In particular, our problem reduces to the study of an ODE in an open subset of \(\R^3\). 
 
 If \(\mu=\mu_{a,b,c}\in \mathcal{O}\), then the norm of \(\mu\) induced from the inner product in \eqref{Eq:InnerProduct} reduces to \[\norm{\mu}^2=a^2+b^2+c^2.\]
 
 It will be useful to have an explicit description of an arbitrary derivation \(D\in \Der(\mu)\) of a bracket \(\mu\in \mathcal{O}\). We can obtain this by differentiating the description of an automorphism of \(\mathfrak{n}_4\) given in \cite{VanThuongMO4DULG}.

\begin{lemma}\label{Derivations}
Let \(\mu =\mu_{a,b,c}\in \mathcal{O}\).  Then, any derivation \(D\in \Der(\mu)\) has the form \[D=\begin{pmatrix}\alpha& 0 & 0 & 0\\ * & \beta &0&0\\ * &a\gamma+b\beta&\alpha+\beta &0 \\ *&*&c\gamma &2\alpha+\beta\end{pmatrix},\qquad \alpha,\beta,\gamma\in \R,\] where each \(*\) is an arbitrary real number.
\end{lemma}

\section{The Bach Flow on a Nilpotent Lie group}\label{BachFlowSection}

Suppose now that \((M^4,g) = (\R^4,g_\mu)\) is a Nilpotent Lie group with left-inavriant metric and that we have a \((\R^4,\cdot_\mu)\)-invariant solution \((g(t))_{t\in(a,b)}\) to the Bach flow with \(g(0)=g_{\mu}\) (a solution \((g(t))_{t\in(a,b)}\) is \((\R^4,\cdot_\mu)\)-invariant if it is for all \(t\in (a,b)\)).  Since the scaler curvature is constant, we find that \(\operatorname{B}(g(t))(0)\) satisfies the ODE \begin{equation}\label{BachFLowODE}
\frac{d}{dt}\m{\cdot}{\cdot}_t = \operatorname{B}(\m{\cdot}{\cdot}_t), \qquad \m{\cdot}{\cdot}_0=g_\mu(0).
\end{equation}
Here, we have written \(\operatorname{B}(\m{\cdot}{\cdot}_t)=\operatorname{B}(g(t))(0)\).  Conversely, given a solution \((\cdot,\cdot)_t\) to the ODE \eqref{BachFLowODE} we obtain a \((\R^4,\cdot_\mu)\)-invariant solution \(g(t)\) to the Bach flow by defining \(g(t)=g_{\mu,(\cdot,\cdot)_t}\) for all \(t\).  By the usual existence and uniqueness of ODE's, we are guaranteed a unique \((\R^4,\cdot_\mu)\)-invariant solution with a maximal interval of existence \cite{LauretGFATSOHS}.  The need for this reasoning is that uniqueness of the Bach flow is an open problem on a general manifold.

Let us fix once and for all initial metric \(g_0 = g_{\mu_0}\) which is invariant under the Nilpotent Lie group \((\R^4,\cdot_{\mu_0})\). Note that we have we have tacitly assumed that the initial inner product on \(T_0\R^n\) is the standard inner product, but this is not an issue since up to isometry this is always the case (see \S\ref{4DNilpotentLieAlgebras}).  

 By the above discussion, there is a unique curve of left-invariant metrics of inner products \(\m{\cdot}{\cdot}_t\) on \(\R^4\) satisfying \eqref{BachFLowODE} which corresponds to the unique \((\R^4,\cdot_{\mu_0})\)-invariant solution of the Bach flow beginning at \(g_{\mu_0}\). It follows that for each \(t\) there is a \(h(t)\in \GL_4(\R)\) such that \[\m{\cdot}{\cdot}_t = \m{h(t)\cdot}{h(t)\cdot}.\]
One can show that the one parameter family of matrices \(h(t)\) can be chosen to be a smooth curve (see \cite[Section 4.1]{LauretGFATSOHS} or Proposition~\ref{Equiv} below).  The corresponding curve of brackets is given by \(\mu(t)= h(t)\cdot \mu_0\).  The bracket flow, introduced by Lauret to study Ricci flow on homogeneous manifolds in \cite{LauretRFOHM}, is motivated by considering what equation the curve \(\mu(t)\in \GL_4(\R)\cdot \mu_0\) satisfies. For more examples of the bracket flow technique see \cite{LauretTRFFSCN,BohmLafuenteIHRF,StanfieldPHCF}.

\begin{definition}
Let \((\R^4,g_{\mu_0})\) as above. The bracket Bach flow is the ODE \begin{equation}\label{BracketFLow}
\frac{d}{dt}\mu = \frac{1}{2}\pi (B_\mu)\mu,\qquad \mu(0)=\mu_0,
\end{equation}
where \(B_\mu\) is defined by \(\operatorname{B}(g_{\mu(t)})(0)=g(0)(B_\mu\cdot,\cdot)\) is the Bach tensor determined by \(\mu\) and \(\pi\) is the representation defined by \[\pi(A)\mu := A\mu - \mu(A\cdot,\cdot)-\mu(\cdot,A\cdot), \qquad\forall A\in \mathfrak{gl}_n.\]
\end{definition}

Note that \(\pi\) is the derivative of the \(\GL_n(\R)\) representation on \(\Lambda^2(\R^n)^*\otimes\R^n\) defined in \eqref{Eq:GLRep}.  Moreover, the solution \(\mu(t)\) of \eqref{BracketFLow} remains in the orbit \(\GL_4(\R)\cdot \mu_0\) since \(\frac{1}{2}\pi(B_{\mu})\in T_{\mu}(\GL_4(\R)\cdot \mu)\) (see \cite[Lemma 3.2]{LauretRFOHM}).
 
Beginning at the inital metric \(g_{\mu_0}\) we now have two families of Riemannian manifolds\begin{equation*}
(\R^4 ,g(t))\qquad (\R^4,g_{\mu(t)})
\end{equation*}
where \(g(t)\) is the unique \((\R^n,\cdot_{\mu_0})\)-invariant solution of the Bach flow and \(\mu(t)\) is the solution of the bracket flow \eqref{BracketFLow}.  The next proposition shows that these are equivalent in a precise way.

\begin{prop}[\cite{LauretTRFFSCN}, Theorem 5.1]\label{Equiv}
Let \((\R^4,g(t))\), \((\R^4, g_{\mu(t)})\) be solutions of the homonegeous Bach flow and the bracket flow respectively.  Then, there exists a family of isomorphisms \(
\varphi(t):(\R^4,\cdot_{\mu_0})\to (\R^4, \cdot_{\mu(t)})\) such that \[g(t)=\varphi(t)^*g_{\mu(t)} \quad\forall t.\]  Moreover, \(h(t)=d\varphi(t)\) satisfies 
\begin{enumerate}
\item \(\m{\cdot}{\cdot}_t=\m{h\cdot}{h\cdot}\)
\item \(\mu(t)= h\mu_0(h^{-1}\cdot,h^{-1}\cdot)\).
\end{enumerate}
\end{prop}
For a proof of Proposition~\ref{Equiv}, one should consult \cite[Theorem 5.1]{LauretTRFFSCN} (note that the proof in \cite{LauretTRFFSCN} is for the Ricci flow however only symmetry of the Ricci tensor is used).
In particular, Proposition~\ref{Equiv} shows that the solutions of \eqref{Eq:BachFLow} and \eqref{BracketFLow} have the same maximal interval of existence and the same curvature (see the Remark after Theorem 3.3 in \cite{LauretRFOHM}).

Recall from \S\ref{4DNilpotentLieAlgebras} that if \(\lambda \in \operatorname{O}(4)\cdot \mu\) then the metrics \(g_\mu\) and \(g_\lambda\) are isometric.  This is due to the fact that if \(h\in \GL_4(\R)\) gives rise to the inner product \((\cdot,\cdot)=\m{h\cdot}{h\cdot}\) and \(k\in \operatorname{O}(4)\) then \(hk\in \GL_4(\R)\) gives rise to the same inner product.  It will be useful to exploit this \(\operatorname{O}(4)\) equivariance when studying the bracket flow \eqref{BracketFLow}.  B\"{o}hm and Lafuente describe this as a refinement of Uhlenbeck's trick of moving frames (see Sections 2 and 3 in \cite{BohmLafuenteIHRF}).  

\begin{prop}[\cite{BohmLafuenteIHRF}, Proposition 3.1]
Let \(R:\GL_4(\R)\cdot \mu_0 \to \mathfrak{so}(4)\) be a smooth map and let \(\mu(t),\bar{\mu}(t)\) denote respectively solutions to the bracket flow \eqref{BracketFLow} and to the modified bracket flow \begin{equation}\label{GuagedBracketFlow}
\frac{d\bar{\mu}}{dt} = \frac{1}{2}\pi(B_\mu - R_\mu)\bar{\mu}, \qquad \bar{\mu}(0)=\mu_0.
\end{equation}
Then, there is a smooth curve \((k(t))\subset \operatorname{O}(4)\) such that \(\bar{\mu} = k\cdot \mu\).
\end{prop}

In particular, the solutions \(\mu,\bar{\mu}\) to the bracket flow \eqref{BracketFLow} and the gauged bracket flow \eqref{GuagedBracketFlow} have the same maximal interval of existence and the same curvature.

We now want to study the bracket flow for \((\R^4,g_\mu)\) when \((\R^4,\mu)\simeq \mathfrak{n}_4\). 

\begin{prop}\label{BachTensor}
Let \(\mu=\mu_{a,b,c}\in \mathcal{O}\).  Then, the endomorphism \(B_\mu:\R^4\to\R^4\) defined by \(\operatorname{B}(g_\mu)(0)=g(B_\mu\cdot,\cdot)\) is given by the following matrix 
\begin{equation}
B_\mu = \begin{pmatrix}
b_{1} & 0 & 0 & 0 \\
0 & b_{2} & b_{5}  & 0 \\
0 & b_{5}  & b_{3} & b_{6}  \\
0 & 0 & b_{6} & b_{4}
\end{pmatrix},
\end{equation}
where 
\begin{align*}
b_{1} &=  \frac{1}{8} \left(4 a^4+ 8a^2 b^2-a^2c ^2+4b^4+8b^2c^2+4c^4\right), \\
b_{2} &=  \frac{1}{24} \left(12 a^4+ 24 a^2b^2-a^2c ^2+12 b^4+8 b^2 c ^2-4 c ^4\right) ,\\
b_{3} &=  \frac{-1}{24} \left(20 a^4-a^2 c ^2+24a^2 b^2+4 b^4-8 b^2 c ^2-12 c ^4\right) , \\
b_{4} &=\frac{1}{24} \left(-4 a^4+3 a^2 c ^2-8 a^2 b^2-20 \left(b^2+c ^2\right)^2\right),\\
b_{5} &= \frac{2}{3} b c  \left(a^2+b^2+c^2\right),  \\
b_{6} &= -\frac{2}{3} a b \left(a^2+b^2+c ^2\right) ,
\end{align*}
\end{prop}

\begin{remark}
From the equations in Proposition~\ref{BachTensor} we can see the following:
\begin{enumerate}
\item We can see explicitly the rescaling formula for the Bach tensor \(B_{c\cdot \mu}=c^4B_\mu\).
\item The expressions for \(b_5\) and \(b_6\) show that \(B_\mu\) is diagonal if and only if \(b=0\).
\item It is interesting to note that \[b_1 = \abs{W_\mu}^2\ge 0.\]  It is not clear to us why this is the case.
\end{enumerate}
\end{remark}

We have noted in \S\ref{4DNilpotentLieAlgebras} that the solution \(\mu(t)\) to the bracket flow \eqref{BracketFLow} beginning at \(\mu_0\in\mathcal{O}\) may not remain in \(\mathcal{O}\).  However, we have also seen in \S\ref{4DNilpotentLieAlgebras} that for any \(\mu \in \mathcal{N}_4\), the orbit \(\operatorname{O}(4)\cdot \mu\) intersects \(\mathcal{O}\).  That is to say, for each \(t\) we can find a \(k(t)\in \operatorname{O}(4)\) such that \(k(t)\cdot \mu(t) \in \mathcal{O}\).  Since \(\mu\) and \(k\cdot\mu\) determine isometric Riemannian manifolds for \(k\in \operatorname{O}(4)\), we may use \(\operatorname{O}(4)\) to gauge our flow, readjusting at each point in time to ensure that the solution remains in \(\mathcal{O}\).  We formalise this in Proposition~\ref{PreservingSymmetries} below by appealing to Proposition~\ref{GuagedBracketFlow} in \S\ref{BachFlowSection}.

\begin{prop}\label{PreservingSymmetries}
Let \(\mu_0=\mu_{a_0,b_0,c_0}\in \mathcal{O}\) and let \(B_{\mu(t)}\) be the Bach endomorphism along the solution \(\mu(t)\) of \eqref{BracketFLow}.  Then, the solution \(\bar{\mu}(t)\) of the gauged bracket flow \eqref{GuagedBracketFlow} remains in \(\mathcal{O}\) where the gauging, \(R_\mu\), is given by \[R_\mu = \begin{pmatrix}
0 & 0 & 0 & 0 \\
0 & 0& b_{5}  & 0 \\
0 & -b_{5}  &0 & b_{6}  \\
0 & 0 & -b_{6} & 0
\end{pmatrix}.\]
\begin{proof}

To show that the symmetries are preserved, it suffices to show that \(\dot{\mu}_{ij}^k=0\) whenever \(i<j\) and \((i,j,k)\notin\{(1,2,3),(1,2,4),(1,3,4)\}\) since then \(\mu_{ij}^k\) will solve the system \(\dot{u}=0,u(0)=0\) and hence be \(u\equiv 0\) by uniqueness.  Here a dot denotes a derivative with respect to time (i.e \(\dot{}:=d/dt)\).

The effect of gauging is that \(B_\mu-R_\mu=L_\mu\) is lower triangular for all \(t\).  With respect to the basis, \(\{e_i\}\), \eqref{GuagedBracketFlow} is \begin{multline}\label{BracketFLowinCoords}
\dot{\mu}_{ij}^k = \frac{d}{dt}\m{\mu(e_i,e_j)}{e_k} = \m{L_\mu \mu(e_i,e_j) - \mu(L_\mu e_i,e_j)-\mu(e_i,L_\mu e_j)}{e_k}\\
=\sum_{l=1}^4\big(\mu_{ij}^lL_l^k - L_i^l\mu_{lj}^k -  L_j^l\mu_{il}^k\big)=\sum_{i<j\le l\le k}\big(\mu_{ij}^lL_l^k - L_i^l\mu_{lj}^k -  L_j^l\mu_{il}^k\big).
\end{multline}
(Note that we are only summing over \(l\)). The last equality follows since \(B_l^k=0\) for \(l>k\) since it is lower triangular and \(\mu_{il}^k=-\mu_{li}^k=0\) for \(l>k\) by our choice of structure constants.  For \(k=1,2\) the right hand side is zero since each term will have a factor of \(\mu_{ij}^1\) or \(\mu_{ij}^2\), all of which are equal to \(0\).  If \(k=3\), then the only triples \((i,j,k)\) we need to check are \((i,3,3)\) for \(i=1,2\).  But \[\dot{\mu}_{i3}^3 = \sum_{i<3\le l\le 3}\big(\mu_{ij}^lL_l^k - L_i^l\mu_{lj}^k -  L_j^l\mu_{il}^k\big) = \mu_{i3}^3L_3^3 - L_i^3\mu_{33}^3 -  L_3^3\mu_{i3}^3=0.\]
If \(k=4\) then we must check \((2,j,4)\) for \(j=3,4\).  This is \[\dot{\mu}_{2j}^4 = \sum_{2<j\le l\le 4}\big(\mu_{ij}^lL_l^k - L_i^l\mu_{lj}^k -  L_j^l\mu_{il}^k\big) =\mu_{2j}^lL_l^4 - L_2^l\mu_{lj}^4 -  L_j^l\mu_{2l}^4 = 0 \] since \(\mu_{2j}^k=0\) for \(j>2\) and \(\mu_{jl}^4=0\) for \(j,l>2\).
\end{proof}
\end{prop}

Therefore, the bracket flow \eqref{BracketFLow} is equivalent to the following ODE for \(\mu=\mu_{a,b,c}\in \mathcal{O}\): 
\begin{align}\label{ODES}
a' &= \frac{-a}{48}(44a^4+72a^2b^2-5a^2c^2+28b^4+24b^2c^4-4c^4)\nonumber \\
b' &= \frac{-b}{48}(60a^4+104a^2b^2+57a^2c^2+44b^4+104b^2c^2+60c^4)\\
c' &= \frac{-c}{48}(-4a^4+24a^2b^2-5a^2c^2+28b^4+72b^2c^2+44c^4)\nonumber
\end{align}

With these in hand we are in a position to study the long time behaviour of the flow.

\begin{lemma}\label{Evolution}
If \(\mu_0\in \mathcal{O}\), the following evolutions hold along the gauged bracket flow:
\begin{align}
 &\frac{d}{dt} \log \frac{a}{c}= (c^2-a^2)\norm{\mu}^2,\\
 &\frac{d}{dt}\frac{b^2}{a^2} \le \frac{-2}{3}\frac{b^2}{a^2} \norm{\mu}^4,\\
 &\frac{d}{dt}\norm{\mu}^2 \le \frac{-1}{12} \norm{\mu}^6.
\end{align}
\begin{proof}
\begin{enumerate}
\item Follows from a direct calculation using \eqref{ODES}\begin{equation*}
\frac{d}{dt}\log \frac{a}{c} = \frac{d}{dt}\log a- \frac{d}{dt}\log c = \frac{a'}{a}-\frac{c'}{c}= (c^2-a^2)\norm{\mu}^2.
\end{equation*}
\item Taking the derivative gives \[\frac{d}{dt}\frac{b^2}{a^2} = 2 \frac{b}{a} \frac{a b' - a'b}{b^2}=2 \frac{b^2}{a^2} \bigg(\frac{b'}{b}-\frac{a'}{a}\bigg). \] 
The claim then follows since \[\bigg(\frac{a'}{a}- \frac{b'}{b}\bigg) = \frac{-1}{24}\big(8\norm{\mu}^4 + 3 c^2( 5a^2+b^2+c^2)  \big)\le \frac{-1}{3}\norm{\mu}^4.\]

\item Similarly to 2, \[48\frac{d}{dt}\norm{\mu}^2 = -44\norm{\mu}^6 + 3a^2 c^2 (47a^2+53b^2+47c^2).\]  Now, we observe that by the multinomial theorem \[\norm{\mu}^6 =(a^2+b^2+c^2)^3 \ge a^6 +c^6 +3a^4 c^2 +3a^2c^4 +6a^2b^2c^2.\]  Hence, \begin{align*}
-40\norm{\mu}^6 + 3a^2 c^2 (47a^2+53b^2+47c^2)&\le -40a^6 +21a^4 c^2 +21a^2c^4 -40c^6  \\ &= (a^2+c^2)(-40a^4+61a^2c^2-40c^4)\\&\le (a^2+c^2)(-40a^4+80a^2c^2-40c^4) = -40(a^2+c^2)(a^2-c^2)^2\le 0.
\end{align*}
Therefore, \[\frac{d}{dt}\norm{\mu}^2 \le \frac{-1}{48}\cdot 4\norm{\mu}^6 = \frac{-1}{12}\norm{\mu}^6.\]
\end{enumerate}
\end{proof}
\end{lemma}
Lemma~\ref{Evolution} gives the following important corollaries.
\begin{corollary}
The Bach flow on a four dimensional simply connected Nilpotent Lie group is immortal, that is, the maximal interval of existence contains \((0,\infty)\).
\begin{proof}
Since \(\norm{\mu}^2\) decreases along the flow, it remains within the closed ball of radius \(\norm{\mu_0}\) which is a compact set.  
\end{proof}
\end{corollary}

\begin{corollary}\label{Decay}
Let \(\mu(t)\) be a solution to \eqref{BracketFLow} in \(\mathcal{N}_4\).  Then,  \[\norm{\mu(t)}^2 \le \frac{\sqrt{6}}{\sqrt{t}}.\]
\begin{proof}
This follows since \(\norm{\mu(t)}^2\) is a sub-solution of \[x' = \frac{-1}{12}x^3, \qquad x(0)=\norm{\mu_0}^2\]  the solution of which is easily seen to be \(x(t) = \sqrt{6}(t+A)^\frac{-1}{2}\le t^{-1/2}\sqrt{6} \).
\end{proof}
\end{corollary}

Corollary~\ref{Decay} shows that the bracket flow beginning at an arbitrary bracket \(\mu_0\in \mathcal{O}\) converges to the Euclidean bracket as \(t\to \infty\).

\section{Self Similar Solutions to the Bach Flow}\label{Sec:Solitons}

In this section we study \emph{solitons} of the Bach flow.  A Bach soliton is a Riemannian manifold \((M^4,g)\) such that  \begin{equation}\label{BachSoliton}
\operatorname{B}(g)+\frac{\Delta s}{12} = \lambda g+\mathcal{L}_X g,
\end{equation}
for a constant \(\lambda \in \R \), and a complete vector field, \(X\in \mathfrak{X}(M)\) (\(\mathcal{L}_X\) is the Lie derivative in the direction \(X\)).  The soliton is called expanding, steady or shrinking if \(\lambda >0\), \(\lambda =0\) or \(\lambda <0\) respectively. If the vector field arises as the gradient of a potential function \(u\in C^\infty(M)\), then \eqref{BachSoliton} becomes \begin{equation}
\operatorname{B}(g) +\frac{\Delta s}{12}= \lambda g+2\operatorname{Hess}u.
\end{equation}
In this case, we say the soliton is a gradient soliton. For homogeneous solitons, the Laplacian term in \eqref{BachSoliton} vanishes since the scalar curvature is constant.  Homogeneous gradient Bach solitons were studied by Griffin in \cite{GriffinGAOSOHM}. Griffin applies Theorem 3.6 from \cite{PetersenWylieROHGSMARE} to reduce the study to Riemannian products of manifolds of dimension less than \(4\). Proposition~\ref{Splitting} below is a simple corollary of \cite[Theorem 3.6]{PetersenWylieROHGSMARE} which we will use to deduce that any non-product solitons we find cannot be gradient.

\begin{prop}[\cite{PetersenWylieROHGSMARE}, \cite{GriffinGAOSOHM}]\label{Splitting}
Let \((M^4,g)\) be a homogeneous Riemannian manifold. If there is a non-constant function \(u\in C^\infty(M)\) such that \[\operatorname{B}(g)=\lambda g+2\operatorname{Hess}u ,\qquad \lambda\in \R,\] then \(M\) splits isometrically as a product \((M,g)=(M'\times \R^k,g'\oplus \overline{g})\) where \(\overline{g}\) is the Euclidean metric and \(u\) is a function on the Euclidean factor.
\end{prop}

Bach solitons correspond to solutions of the Bach flow which evolve by scaling and pullback by diffeomorphisms \cite[Theorem 4.10]{LauretGFATSOHS}.

If \((M^4,g)=(\R^n,g_\mu)\) is a simply connected Nilpotent Lie group, then an analogous condition to the metric evolving self similarly under the Bach flow is that the solution \(\mu(t)\) to the bracket flow evolves only by scaling (see the discussion before Theorem 6 in \cite{LauretGFATSOHS}).  If this is the case, the bracket \(\mu=\mu(0)\) satisfies \begin{equation}\label{AlgSoliton}
B_\mu =\lambda I + D, \qquad \lambda\in \R, D\in \Der(\mu).
\end{equation}

A simply connected Nilpotent Lie group \((\R^n,g_\mu)\) whose Bach endomorphism satisfies \eqref{AlgSoliton} is called an \emph{algebraic Bach soliton}.  Algebraic solitons were introduced by Lauret to study Ricci Nilsolitons (Ricci solitons on Nilpotent Lie groups) in \cite{LauretRSHN}. In \cite{LauretGFATSOHS} Lauret generalises the the bracket flow technique and algebraic solitons to a large class of geometric structures. An important observation made by Lauret is that an algebraic soliton is indeed a soliton in the sense of \eqref{BachSoliton}.  For convenience, we have summarised this for the Bach tensor in Proposition~\ref{prop:LauretASol} below.

\begin{prop}[Theorem 6, \cite{LauretGFATSOHS}]\label{prop:LauretASol}
For a simply connected Nilpotent Lie group \((\R^n,g_{\mu_0})\) the following are equivalent:
\begin{enumerate}
\item  The solution to the bracket flow \eqref{BracketFLow} beginning at \(\mu_0\) is given by \[\mu(t) = \lambda(t)\cdot \mu_0, \quad \text{for some }\lambda(t)>0,\, \lambda(0)=1. \]
\item The operator \(B_\mu\) associated to the Bach tensor satisfies \eqref{AlgSoliton}.
\end{enumerate}
Moreover, whenever either of these conditions hold, the Riemannian manifold \((\R^n,g_{\mu_0})\) is a Bach soliton.
\end{prop} 
\begin{remark}
If \(\lambda<0\) in \eqref{AlgSoliton} then the corresponding bracket flow from Proposition~\ref{prop:LauretASol} will shrink homothetically towards \(0\).  Since \(g_{\lambda\cdot\mu}=\lambda^{-2}g_\mu\) this corresponds to an expanding Bach soliton.  
\end{remark}

Note that in general, not all solitons are algebraic solitons.  Since we have an explicit description of what a derivation \(D\in \Der(\mu)\) looks like for a bracket \(\mu\in \mathcal{O}\), the notion of an algebraic soliton reduces our search for a soliton to simply solving a system of equations in terms of \(a,b,c,\) and the components of \(D\) (note that since \(B_\mu\) is trace free, \(\lambda=-\tr D/4\)).  In fact, we can reduce the difficulty of this system further.

\begin{lemma}
If \(\mu=\mu_{a,b,c}\in \mathcal{O}\) is an algebraic soliton, then \(B_\mu\) is diagonal.
\begin{proof}
Since \(D=B_\mu - \lambda I\), the derivation \(D\in \Der(\mu)\) must be symmetric as the difference of two symmetric matrices.  But by Lemma~\ref{Derivations}, \(D\) is lower triangular.  Hence, \(D\) must be diagonal and so \(B_\mu=\lambda I+D\) must also be diagonal.
\end{proof}
\end{lemma}

Since \(B_\mu\) is diagonal if and only if \(b=0\), this allows us to set \(b=0\) when searching for algebraic solitons.  By Lemma~\ref{Derivations}, a diagonal derivation has eigenvalues \(\alpha,\beta, \alpha+\beta,2\alpha+\beta\) for some \(\alpha,\beta\in \R\).  Therefore, existence and uniqueness of algebraic solitons reduces to existence and uniqueness of solutions to a set of polynomials in \(a,b,\alpha,\beta\).

\begin{theorem}\label{Thm:BachSoliton}
The bracket \(\mu\in \mathcal{O}\subset \mathcal{N}_4\) given by \[\mu(e_1,e_2)=e_3 \quad \mu(e_1,e_3)=e_4\] is a Bach soliton.  Moreover, this soliton is a non-gradient expanding soliton and is the unique algebraic Bach soliton up to isometry and scaling within the orbit \(\GL_4(\R)\cdot \mu \subset \mathcal{N}_4\).
\begin{proof}
 By setting \(a=c=1\) and \(b=0\) in Proposition~\ref{BachTensor} it is not difficult to check that \(\alpha=-7/12\) and \(\beta=-7/6\) solves \eqref{AlgSoliton}. This gives \(\lambda = -21/16 <0\) so the soliton is expanding.  The soliton is not of gradient type due to Proposition~\ref{Splitting}.
 
To see uniqueness we assume that \eqref{AlgSoliton} holds for \(a,c,\alpha,\beta\in \R\) such that \(\mu_{a,0,c}\in\mathcal{O}\) and \(a\neq c\). We then show that this leads to a contradiction.
\end{proof}
\end{theorem}

\begin{remark}
Up to isometry any left-invariant metric on \(\R\times H^3\) can be described by a bracket \(\mu=\mu_{a,b,c}\) with \(b=c=0\) and \(a>0\).  By Proposition~\ref{prop:LauretASol} these must be solitons since they trivially evolve by scaling.  Moreover, by \cite[Proposition 4.14]{GriffinGAOSOHM} these are also non-gradient.
\end{remark}

\section{Normalised Bach Flow and Convergence}\label{Sec:NormalisedFLows}

Similarly to the Ricci flow case \cite[Section 7]{LauretTRFFSCN}, one may consider a normalised Bach flow:
 \begin{equation}\label{rNormBachFlow}
\frac{\partial}{\partial t}g(t) = \operatorname{B}(g(t)) +r(t)g(t), \quad g(0)=g_0,
\end{equation}
for some normalisation function \(r:[0,T)\to\R\).

In the varying brackets perspective the analogue of \eqref{rNormBachFlow} is the following \(r-\)normalised bracket flow equation.

\begin{definition}
An \(r\)-normalised bracket Bach flow for a normalisation function \(r:[0,T)\to \R\) is a curve \((\mu^r(t))\subset \mathcal{N}_4\) such that 
\begin{equation}\label{rNormBracketFlow}
\frac{d}{dt}\mu^r = \frac{1}{2}\pi(B_{\mu^r} -rI){\mu^r},\quad \mu(0)=\mu_0.
\end{equation}
\end{definition}
\begin{remark}
There is an equivalence between \eqref{rNormBachFlow} and \eqref{rNormBracketFlow} which is analogous to the case of un-normalised flows given in Proposition~\ref{Equiv} \cite[Section 7]{LauretTRFFSCN}.
\end{remark}

The usefulness of \eqref{rNormBachFlow} and \eqref{rNormBracketFlow} is that the addition of the \(r\) term allows us to keep a geometric quantity fixed along the flow.  Moreover, the next proposition shows that solutions of the normalised flows only differ from the solutions of the original flows by a scaling and reparametrisation of time.

\begin{prop}\label{Prop:Normalised}
Let \(\mu(t)\) and \(\mu^r(t)\) be solutions of the bracket flow \eqref{BracketFLow} and the \(r-\)normalised bracket flow \eqref{rNormBracketFlow} respectively.  Then, there are functions \(\tau:[0,T)\to [0,T)\), \(\lambda:[0,T)\to \R\) such that \begin{equation}
\mu^r(t) = \lambda(t)\mu(\tau(t)) \qquad \forall t\in [0,T).
\end{equation}
The fucntions \(\tau\) and \(\lambda\) are the solutions of the ODE's \begin{equation}\label{RescalingForNormalisedFlow}
\tau ' = \lambda^4 ,\quad \tau(0)=0 \qquad \qquad \lambda' = \frac{1}{2}r\lambda ,\quad \lambda(0)=1.
\end{equation}
\begin{proof}
Let \(\mu(t)\) be a solution to \eqref{BracketFLow} and define \[\mu^r(t)=\lambda(t)\mu(\tau(t))\] where \(\tau,\lambda\) are the solutions of \eqref{RescalingForNormalisedFlow}.  Clearly \(\mu^r(0) = \mu_0\).  Differentiating gives \begin{multline*}
\frac{d}{dt}\mu^r(t)=\frac{d}{dt}(\lambda(t)\mu(\tau(t))) = \lambda'\mu(\tau(t))+\lambda \tau'\frac{d}{dt}\bigg|_{\tau(t)}\mu=\frac{1}{2}r\lambda\mu(\tau(t)) +\frac{1}{2}\lambda^5\pi(B_\mu)\mu \\
= \frac{1}{2}r\mu^r +\frac{1}{2} \pi(B_{\mu^r})\mu^r = \frac{1}{2}\pi(B_{\mu^r} - rI)\mu^r
\end{multline*}
where we have used that the Bach tensor scales by \(B_{\lambda\cdot\mu} = \lambda^4B_\mu.\)
\end{proof}
\end{prop}

\begin{corollary}\label{ScaleInvariant}
Let \(u:\Lambda^2(\R^4)^*\otimes \R^4\to \R\) be scale invariant (i.e. \(u(\lambda\mu)=u(\mu)\) for any \(\lambda\in \R^*\)) and smooth away from \(0\).  Then \(u\) increases (resp. decreases) along a solution of the bracket flow if and only if it \(u\) increases (resp. decreases) along a solution if the normalised bracket flow.
\begin{proof}
Let \(\mu,\mu^r\) be solutions to the bracket flow and \(r\)-normalised bracket flow respectively.  If \(\tilde{u},u^r\) denote the restrictions of \(u\) to \(\mu,\mu^r\) then Proposition~\ref{Prop:Normalised} and scale invariance of \(u\) imply that \(u^r(t)=\tilde{u}(\tau(t))\).  Since \(\tau\) is increasing, the claim follows.
\end{proof}
\end{corollary}

In particular, if \(\norm{\mu_0}=2\) and  \[r =-\frac{1}{4}\m{\pi(B_\mu)\mu}{\mu},\]  then \(\norm{\mu(t)}\equiv 2\) since \[\frac{d}{dt}\norm{\mu}^2 = 2 \m{\frac{1}{2}\pi(B_\mu-rI)\mu }{\mu} = \m{\pi(B_\mu)\mu}{\mu} - \frac{ \m{\pi(B_\mu)\mu}{\mu}\norm{\mu}^2 }{\norm{\mu}^2}=0. \]
The scalar curvature is also constant since \(s_\mu = -\norm{\mu}^2/4\) on Nilpotent Lie groups (c.f. 7.39 in \cite{BesseEM}).

Since the quantities \(a/c\) and \(b^2/a^2\) are scale invariant, we can use Lemma~\ref{Evolution} and Corollary~\ref{ScaleInvariant} to determine the behaviour of the normalised flow.

\begin{theorem}\label{Thm:Convergence}
Let \(\mu_0\in \mathcal{O}\) with \(\norm{\mu_0}=2\).  If \(\mu(t)\) is the solution to the \(r\)-normalised Bach flow \eqref{rNormBracketFlow} for \(r=- \frac{1}{4}\m{\pi(B_\mu)\mu}{\mu}\) then \(g_{\mu(t)}\to g_{\mu_\infty}\) uniformly on compact subsets of \(\R^4\) as \(t\to \infty\) where \(\mu_\infty\) is (up to scaling) the algebraic Bach soliton given in Theorem~\ref{Thm:BachSoliton}.
\begin{proof}
Let \(\mu_0=\mu_{a_0,b_0,c_0}\in \mathcal{O}\) with \(a_0^2+b_0^2+c_0^2=4\).  We show that \(\mu(t)\to \mu_\infty\), \cite[Proposition 2.1]{LauretTRFFSCN} then implies \(g_{\mu(t)}\to g_{\mu_\infty}\) uniformly on compact subsets of \(\R^4\) as \(t\to \infty\). Showing \(\mu(t)\to \mu_\infty\) amounts to showing that \(a,c\to 2\) and \(b\to 0\) as \(t\to \infty\).

We first show that \(b\to 0\) as \(t\to \infty\).  By Lemma~\ref{Evolution} and Corollary~\ref{ScaleInvariant}, we know that \[\frac{d}{dt}\frac{b^2}{a^2}\le \frac{-2}{3}\frac{b^2}{a^2}\tau'<0.\]   Therefore, \(b^2/a^2\) is monotone decreasing and is bounded below by \(0\) so it must converge.  But this implies the derivative of \(b^2/a^2\) must apporach \(0\), which implies \(b^2/a^2\) converges to \(0\) since \(\tau'=\lambda^4\ge 1\) for \(t\) sufficiently large.  Since \(a\le 2\) this implies \(b\to 0\) as \(t\to \infty\).  

Next, we have that \[\frac{d}{dt}\log\frac{a}{c} = \tau'(c^4-a^4)\norm{\tilde{\mu}(\tau(t))} =2\lambda^3(c^4-a^4).\] Here \(\tilde{\mu}\) denotes the un-normalised flow and we have used that \(2=\norm{\mu}=\lambda\norm{\tilde{\mu}}\). Observe that if \(a(t_0)=c(t_0)\) for some \(t_0\ge 0\), then \(a\equiv c\) for all \(t>t_0\) by uniqueness of ODE solutions.  Therefore, if \(a(0)>c(0)\) (resp. \(a(0)<c(0)\)) then we can assume \(a(t)>c(t)\) for all \(t\) (resp. \(a(t)<c(t)\)).  In this case, \(\log\frac{a}{c}\) is a monotone and bounded, and hence must be convergent.  But then \(d/dt \log \frac{a}{c}\to 0\) as \(t\to \infty\) so \(c^4-a^4\to 0\).  Since \(a^2+c^2\to 1\) as \(t\to \infty\), it must hold that \[\lim_{t\to \infty}a =\lim_{t\to \infty} c = 2.\]
\end{proof}
\end{theorem}

Note that Theorem~\ref{Thm:Convergence} implies that the Bach flow \(g(t)=\varphi(t)^*g_{\mu(t)}\) converges in the pointed Cheeger-Gromov sense to \(g_\infty:=g_{\mu_\infty}\) which is a non-trivial, non-gradient algebraic Bach soliton by Theorem~\ref{Thm:BachSoliton}.

\begin{acknowledgements}
I am grateful to Dr Ramiro Lafuente for his support and encouragement throughout this research.
\end{acknowledgements}

\bibliography{BachFlowOnSimplyConnectedNilmanifolds}

\end{document}